\documentclass{article}
\usepackage{amsthm}
\usepackage{graphicx}
\newtheorem{lem}{Lemma}
\newtheorem{thm}[lem]{Theorem} 
\newtheorem{df}[lem]{Definition}

\newtheorem{exa}[lem]{Example} 
 
\newcommand{\ex}{\mbox{ex}}
\title{A characterization of the edge-shelling convex geometries of trees}
\author{Kenji Kashiwabara and Masataka Nakamura\\
Department of Systems Science, The University of Tokyo}
\date{}
\begin{document}

\maketitle

\begin{abstract}
We investigate the class of the edge-shelling convex geometries of trees. The edge-shelling convex geometry of a tree is the convex geometry consisting of the sets of edges of the subtrees. For the edge-shelling convex geometry of a tree, the size of the stem of any rooted circuit is two. The class of the edge-shelling convex geometry of a tree is closed under trace operation.
We characterize the class of the edge-shelling convex geometry of a tree in terms of trace-minimal forbidden minors. 
Moreover, the trace-minimal forbidden minors are specified for the class of convex geometries such that the size of any stem is two.
\end{abstract}

\section{Introduction}

A convex geometry is a combinatorial abstraction of convexity whereas a matroid is an abstraction of linear dependency. Convex geometries are derived from a variety of combinatorial objects, such as posets, affine point configurations, chordal graphs, semi-lattices, and so on \cite{EJ}. The complement of a convex geometry is known as an antimatroid. We investigate convex geometries whose rooted circuits have stems of size 2.

A tree is a connected graph containing no cycle.
 A tree provides two types of convex geometries, that is, the edge-shelling convex geometry of the tree and the vertex-shelling convex geometry of the tree.
The edge-shelling convex geometry of a tree is the convex geometry consisting of the sets of edges of the subtrees.
A subtree of a tree is a connected subgraph.
The vertex-shelling convex geometry of a tree is the convex geometry consisting of the sets of vertices of the subtrees.
The size of any stem of the edge-shelling convex geometry of a tree is two. Every trace of the edge-shelling convex geometry of a tree is also the edge-shelling convex geometry of some tree. On the other hand, the size of some stem of some vertex-shelling convex geometry is not two. Some trace of some vertex-shelling convex geometry is not a vertex-shelling convex geometry. We compare an edge-shelling convex geometry and a vertex-shelling convex geometry in Subsection \ref{subsec:edgevertex}.

In this paper, we characterize the class of the edge-shelling convex geometries of trees. We specify all the trace-minimal forbidden minors. This is our main result, presented in Section \ref{sec:main}. 

A simplicial-shelling convex geometry is the convex geometry arising from a chordal graph.
Any rooted circuit of a simplicial-shelling convex geometry has stems of size 2.
We show an example of a simplicial-shelling convex geometry whose trace does not belong to the class of a simplicial-shelling convex geometry. This will be presented in Section \ref{sec:simplicial}.

Moreover, the trace-minimal forbidden minors are determined for the class of convex geometries such that the size of any stem is two. This class is characterized in the families of rooted sets whereas the class of an edge-shelling convex geometry is characterized in the convex geometries.
This will be presented in Section \ref{sec:eminor}.

Although there are a rich number of forbidden-minor type characterization theorems in matroid theory, there are a few characterization theorems in the theory of convex geometries. Nakamura \cite{N} gives a forbidden minor characterization for the class of convex geometries of the graph search with respect to the operations of deletion and contraction. Okamoto and Nakamura \cite{ON} gives that of the convex geometries of the graph line search. 

\section{Preliminaries}

\subsection{Closure systems}

Let $E$ be a finite ground set throughout this papar.

Since a convex geometry is a special case of closure systems, we begin with defining the closure system.
The closure operator and the extreme operator are defined for a closure system. These operators play an important role for convex geometries.

\begin{df}
A family $\cal K$ of sets on $E$ is a closure system if the following conditions hold.

\begin{enumerate}
\item $E\in {\cal K}$.
\item $A,B\in {\cal K}$ implies $A\cap B \in {\cal K}$.
\end{enumerate}
\end{df}

An element of a closure system is called a {\it closed set}.

For a closure system, define $\tau(X)=\bigcap\{A \in {\cal K}|X\subseteq A\}$ for $X\subseteq E$.
Since $\tau(X)$ is a closed set, it is the smallest set in all the closed sets including $X$. 
$\tau:2^E\to 2^E$ is called the {\it closure operator} of the closure system.

Define $ex(X) = \{x\in X|x \notin \tau(X-x)\}$ for $X\subseteq E$. $\ex:2^E\to 2^E$ is called the {\it extreme operator} of the closure system.

%

For $e\in V$ and $X\subseteq E$, $(X,e)$ with $e\notin X$ is called a {\it rooted set}.
$e$ is called the {\it root} of the rooted set, and $X$ is called the {\it stem} of the rooted set. 
In the literature, $(X\cup\{e\},e)$, other than $(X,e)$, is often called a rooted set. But when we write a rooted set $(X,e)$, we assume that $X$ does not include $e$.

For a closure system ${\cal K}$ on $E$, a rooted set $(X,e)$ is called a {\it rooted circuit} if 
$X$ is a minimal set satisfying $e\in \tau(X)$.

The next lemma gives the relation between the extreme operator and the rooted circuits.

\begin{lem}\label{lem:exstem}
For a closure system, $e\in ex(X)$ if and only if there exists no rooted circuit $(A,e)$ with $A\subseteq X$.
\end{lem}

For a family $\cal C$ of rooted circuits and $T\subseteq E$, the {\it trace}  ${\cal C}:T$ from $\cal C$ on $T$ is defined as $\{(X,e)\in {\cal C}|X \subseteq T, e\in T\}$.

While trace operation is usually defined for a family of sets, we define trace operation for a family of rooted circuits. We will see in Lemma \ref{lem:tracerooted} that both definitions coincide with each other for convex geometries. This definition plays an important role in Section \ref{sec:eminor}.

\subsection{Convex geometries}

A convex geometry is a special case of a closure system.
We can define that a closure system is a convex geometry in various ways.

\begin{df}
A closure system ${\cal K}$ is a {\it convex geometry} if, for any non-empty closed set $X\in {\cal K}$, $\ex(X)$ is non-empty.
\end{df}

An element of a convex geometry is called a {\it convex set}.

It is known that a convex geometry is also defined in terms of the axiom of rooted circuits.

\begin{lem}\cite{D}\label{lem:rootedaxiom}
A family $\cal C$ of rooted sets becomes the set of rooted circuits of a convex geometry if and only if the following two conditions hold.

(1) If $(X,e)$ and $(Y,e)$ belong to $\cal C$ with $X\subseteq Y$, then $X=Y$.

(2) If $(X,e)$ and $(Y,f)$ belong to $\cal C$ with $e\in Y$, then there exists $Z\subseteq X\cup Y-\{e\}$ with $(Z,f) \in \cal C$.
\end{lem}

The next lemma follows from the definitions of trace and rooted circuits.

\begin{lem}\label{lem:tracerooted}
The trace ${\cal C}:T$ from the rooted circuits $\cal C$ of a convex geometry $\cal K$ 
on $T\subseteq E$ is the rooted circuits of a convex geometry $\{T \cap X|X \in {\cal K}\}$.
Note that $T$ may not be a convex set.
\end{lem}

So the class of a convex geometry is closed under trace operation.

\subsection{Convex geometries of stem size 2}

The class of edge-shelling convex geometries of trees is a subclass of convex geometries with stems of size 2.
In this subsection, we consider convex geometries with stems of size 2.
We call such a convex geometry a {\it convex geometry of stem size 2}.
There are many important classes of convex geometries which belong to this class, e.g., the double-shelling convex geometry of a poset, the simplicial-shelling convex geometry of a chordal graph, the edge-shelling convex geometry of a tree, the vertex-shelling convex geometry of a tree, and so on.

The existence of a rooted circuit $(\{x,y\},z)$ may be interpreted as that $z$ is between $x$ and $y$ generally.

We consider a family of rooted circuits $(\{x,y\},z)$ and we discuss whether it satisfies the axiom of the convex geometry or not.
At that time, the trace of such a family on $X\subseteq E$ is defined as the collection of $(\{x,y\},z)$ with $x,y,z\in X$.

The next lemma follows from the case where the size of stem is 2 in Lemma \ref{lem:rootedaxiom}.

\begin{lem}\label{lem:t}
For any convex geometry with the stem size 2, when $(\{x,y\},z)\in {\cal C}$ and $(\{v,z\},u)\in {\cal C}$, then $(\{x,v\},u)\in {\cal C}$, $(\{y,v\},u)\in {\cal C}$, or $(\{x,y\},u)\in {\cal C}$ holds.
\end{lem}

\begin{proof}
By Lemma \ref{lem:rootedaxiom}, there exists a rooted circuit $(X,u) \in {\cal C}$ such that $X\subseteq \{x,y,v\}$. Since the size of the stem is 2, $X=\{x,v\}, \{y,v\},$ or $\{x,y\}$.
\end{proof}

The following lemmas are special cases of Lemma \ref{lem:t}.

\begin{lem}\label{lem:p}
For any convex geometry with the stem size 2, when $(\{x,y\},z)\in {\cal C}$ and $(\{x,z\},u)\in {\cal C}$, then $(\{x,y\},u)\in {\cal C}$ holds.
\end{lem}

\begin{lem}\label{lem:i4}
For any convex geometry with the stem size 2, when $(\{x,y\},z)\in {\cal C}$ and $(\{u,z\},y)\in {\cal C}$, $(\{x,u\},y)\in {\cal C}$ and $(\{x,u\},z)\in {\cal C}$ hold.
\end{lem}

\begin{figure}[ht]
\begin{center}
\begin{minipage}[ht]{4cm}
\begin{center}
\includegraphics[bb=0 0 40mm 30mm]{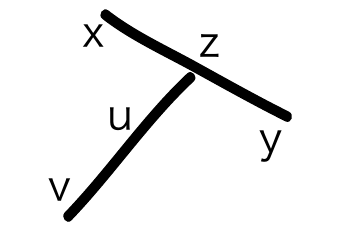}
\caption{Lemma \ref{lem:t}}\label{fig:rule-t}
\end{center}
\end{minipage}
\begin{minipage}[ht]{4cm}
\begin{center}
\includegraphics[bb=0 0 40mm 30mm]{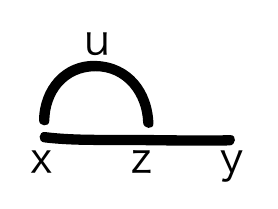}
\caption{Lemma \ref{lem:p}}\label{fig:rule-p}
\end{center}
\end{minipage}
\begin{minipage}[ht]{4cm}
\begin{center}
\includegraphics[bb=0 0 40mm 30mm]{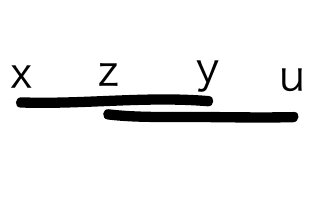}
\caption{Lemma \ref{lem:i4}}\label{fig:rule-i}
\end{center}
\end{minipage}
\end{center}
\end{figure}

In this paper, we illustrate a rooted circuit as a line connecting three points as in Figures \ref{fig:rule-t}, \ref{fig:rule-p}, and \ref{fig:rule-i}. Since these figures do not depict Hasse diagrams, the vertical position of points in these figures has no meaning. When a figure depicts the Hasse diagram of a poset, we write explicitly `the Hasse diagram' in its caption to avoid confusion.

We consider the ordering arising from a convex geometry of stem size 2 and an extreme element as follows.

Fix an element $x\in \ex(E)$ since $\ex(E)$ is not empty.

When $x$ is contained in some rooted circuit, $x$ is contained in the stem of the rooted circuit since $x\in \ex(E)$ by Lemma \ref{lem:exstem}.

For $y,z \in E-\{x\}$, define a binary relation $y>_{x} z$ when $(\{x,z\},y)$ is a rooted circuit. 
Define $x >_{x} y$ for any $y\in E-\{x\}$.

\begin{lem}\label{lem:strict}
For a convex geometry of stem size 2 and $x \in \ex(E)$, binary relation $>_x$ is a strict partial order.
\end{lem}

\begin{proof}
It suffices to show asymmetry and transitivity.

First we show asymmetry. Since the case where $x$ is involved in is easy, we prove the asymmetry which does not involve $x$. Suppose that $y>z$ and $z>y$. Then we have $(\{x,z\},y)\in {\cal C}$ and $(\{x,y\},z)\in {\cal C}$. This contradicts that the trace on $\{x,y,z\}$ is a convex geometry.

Next, we show transitivity. Since the case where $x$ is involved in is easy, we prove the transitivity which does not involve $x$. Suppose that $y>z$ and $z>u$. Then we have $(\{x,z\},y)\in {\cal C}$ and $(\{x,u\},z)\in {\cal C}$. By Lemma \ref{lem:p}, we have $(\{x,u\},y)\in {\cal C}$ and $y>u$.
\end{proof}

\subsection{Edge-shelling convex geometry of a tree}\label{subsec:edgevertex}

The convex geometry arising from the edge shelling of a tree is a convex geometry defined on the edges of the tree. Its convex sets consist of a subset of the edges of the given tree that forms a tree. 
The class of convex geometries arising from the edge shelling of a tree is closed under trace operation. The size of any stem of the edge-shelling convex geometry of a tree is two.

\begin{lem}
For the edge-shelling convex geometry of a tree with rooted circuits ${\cal C}$, $(\{x,y\},z)\in {\cal C}$ if and only if $z$ is on the path between $x$ and $y$ on the tree.
\end{lem}

\begin{proof}
This lemma follows from that $\tau(\{x,y\})$ consists of the edges on the path between $x$ and $y$.
\end{proof}

The vertex-shelling convex geometry of a tree is the convex geometry consisting of the sets of vertices of the subtrees.
The size of any stem of the edge-shelling convex geometry of a tree is two. Every trace of the edge-shelling convex geometry of a tree is also the edge-shelling convex geometry of some tree. On the other hand, the size of some stem of some vertex-shelling convex geometry is not two. Some trace of some vertex-shelling convex geometry is not a vertex-shelling convex geometry.

\begin{exa}
We give an example of the edge-shelling convex geometry and the vertex-shelling convex geometry of a tree. The tree shown in Figure \ref{fig:treeshelling}(left) gives arise to the edge-shelling convex geometry shown in Figure \ref{fig:treeshelling}(middle) and the vertex-shelling convex geometry shown in Figure \ref{fig:treeshelling}(right). A curve connecting three points in these figures depicts a rooted circuit of stem size 2. For example, edge b is on the path between edge a and edge c in the tree. So the edge-shelling convex geometry has rooted circuits $(\{a,c\},b)$ and $(\{a,d\},b)$. The vertex-shelling convex geometry has rooted circuits 
$$ (\{A,C\},B), (\{A,D\},C), (\{A,D\},B), (\{B,D\},C), (\{A,E\},C), (\{A,E\},B), (\{B,E\},C).$$

\begin{figure}[ht]
\begin{center}
\includegraphics[bb=0 0 100mm 45mm]{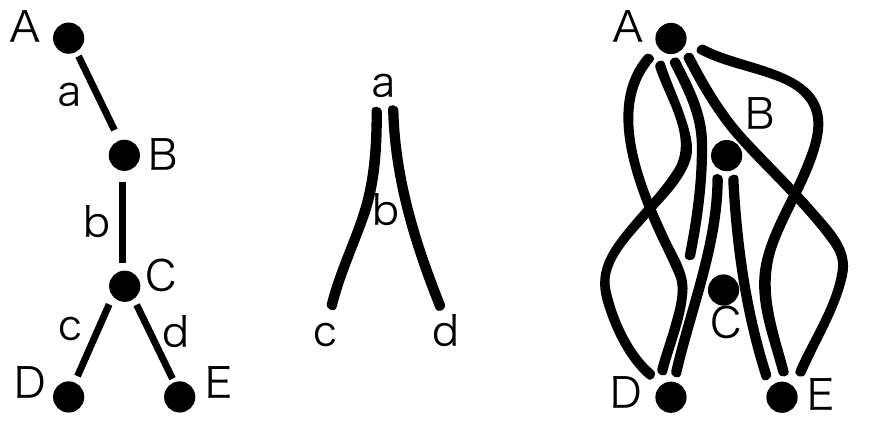}
\caption{Tree, edge-shelling convex geometry, and vertex-shelling convex geometry}\label{fig:treeshelling}
\end{center}
\end{figure}

\end{exa}

The next lemma means that the class of edge-shelling convex geometries is closed under trace operation.

\begin{lem}
The trace of the edge-shelling convex geometry of a tree is also the edge-shelling convex geometry of some tree.
The trace on $A\subseteq E$ corresponds to the contraction of the edges in $E-A$ from the graph.
\end{lem}

\begin{proof}
It is because the contracting the graph does not affect whether $z$ is on the path between $x$ and $y$.
\end{proof}

\section{Characterization of edge-shelling convex geometries}\label{sec:main}

\subsection{Main theorem}

In this section, we present trace-minimal forbidden minors for the class of the edge-shelling convex geometry of a tree. The next theorem is the main theorem of this paper, which gives a forbidden-minor-type characterization for edge-shelling convex geometries.

\begin{thm}\label{thm:main}
For a convex geometry with the rooted circuits ${\cal C}$, it is the edge-shelling convex geometry of a tree if and only if any trace of it coincides with none of following convex geometries.

\begin{itemize}
\item Type i:$(\{y,z\},u)$ on $\{x,y,z,u\}$.
\item Type Y:$(\{x,y\},u),(\{y,z\},u),(\{x,z\},u)$.
\item Type O:$(\{x,z\},y),(\{x,z\},u)$.
\item Type OC:$(\{x,z\},y),(\{x,z\},u),(\{u,z\},y)$.
\item Type B:$(\{x,z\},y),(\{x,y\},u),(\{y,z\},u),(\{x,z\},u)$.
\end{itemize}

\begin{figure}[ht]
\begin{center}
\includegraphics[bb=0 0 140mm 46mm]{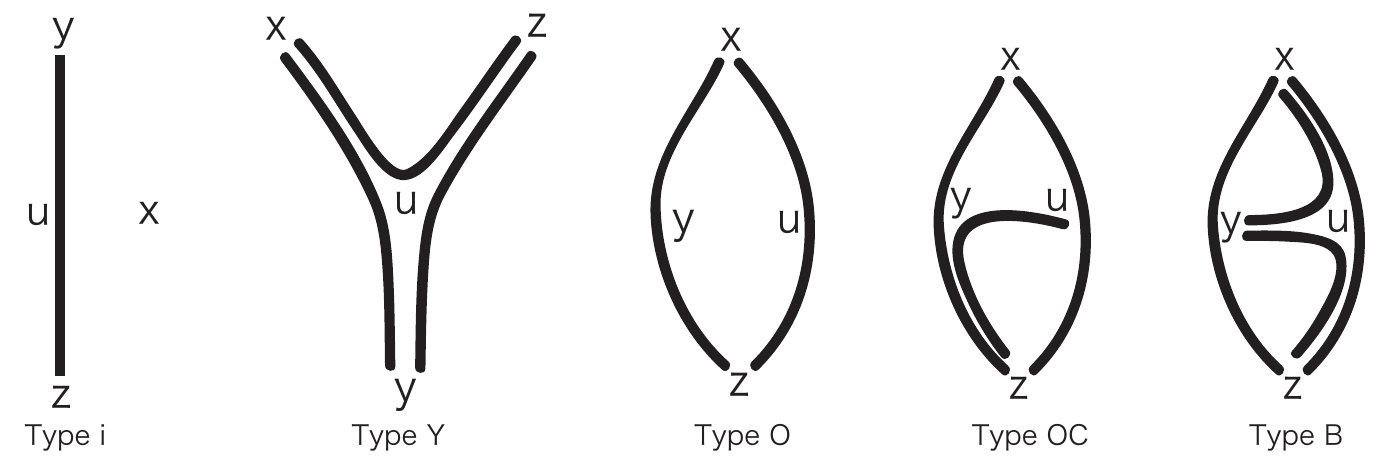}
\caption{Trace-minimal non-edge-shelling convex geometries}\label{fig:forbiddenedge}
\end{center}
\end{figure}
\end{thm}

We call each of convex geometries Type i, Type Y, Type O, Type OC, and Type  B a {\it trace-minimal non-edge-shelling convex geometry}. In fact, each of them is not an edge-shelling convex geometry but any proper minor is an edge-shelling convex geometry.

\begin{lem}
For a convex geometry of stem size 2, it does not have any trace-minimal non-edge-shelling convex geometry if and only if the following rule is satisfied.

Rule Y: For distinct four elements $x,y,z,u$, $(\{y,z\},u)\in {\cal C}$ implies that exactly one of 
$(\{x,y\},u)\in {\cal C}$ and $(\{x,z\},u)\in {\cal C}$ holds.
\end{lem}

\begin{proof}
Since no trace-minimal non-edge-shelling convex geometry in Figure ¥ref{fig:forbiddenedge} satisfies Rule Y, Rule Y is not satisfied when the convex geometry has a trace-minimal non-edge-shelling convex geometry.

Conversely, suppose that the convex geometry does not satisfy Rule Y. Then there exist some four elements which satisfy the assumption of Rule Y and does not satisfy the conclusion of Rule Y. There exist two cases, that is, both satisfy the conclusion and both do not satisfy the conclusion. 

We first consider the case where both conditions satisfy the conclusion. In this case, we may assume that the convex geometry has at least $(\{x,y\},u), (\{y,z\},u), $ or $(\{x,z\},u)$. It may have other rooted circuits in $\{x,y,z,u\}$. For a convex geometry each three set included in $\{x,y,z,u\}$ cannot allow more than one rooted circuit. So if it has an additional rooted circuit on $\{x,y,z,u\}$, the ground set of the rooted circuit must be $\{x,y,z\}$. In this case, the convex geometry obtained by adding each of $(\{x,y\},u),(\{y,u\},x)$ and $(\{u,x\},y)$ is isomorphic to Type B. So the convex geometry must have Type Y or Type B as a trace minor.

Next we consider the case where neither of the two conditions satisfies the conclusion.
The convex geometry has at least $(\{y,z\},u)$. It may have other rooted circuits in $\{x,y,z,u\}$. When $(\{x,u\},y)$ is a rooted circuit, $(\{x,z\},u)$ is a rooted circuit by Lemma \ref{lem:i4}, which contradicts the assumption of this case. When $(\{x,z\},u)$ is a rooted circuit, $(\{x,z\},u)$ is a rooted circuit by Lemma \ref{lem:p}, which contradicts the assumption of this case. When $(\{y,u\},x)$ is a rooted circuit, $(\{y,z\},x)$ should be a rooted circuit by Lemma \ref{lem:p}. So the convex geometry has Type i, Type O, Type OC, or Type B.
\end{proof}

\begin{lem}
For a convex geometry of stem size 2, Rule Y implies Rule O.

Rule O: When $(\{x,z\},y)\in {\cal C}$ and $(\{x,z\},u)\in {\cal C}$ hold, $(\{x,y\},u)\in {\cal C}$ or $(\{x,u\},y)\in {\cal C}$ holds.
\end{lem}

\begin{proof}
Suppose that $(\{x,y\},u)\notin {\cal C}$ and $(\{x,u\},y)\notin {\cal C}$. 
By applying Rule Y to $y$ and $(\{x,z\},u)$, we have $(\{x,z\},u)\in {\cal C}$.
By applying Rule Y to $u$ and $(\{x,z\},y)$, we have $(\{x,y\},u)\in {\cal C}$.  By $(\{x,z\},u)\in {\cal C}$ and $(\{x,y\},u)\in {\cal C}$, the trace on $\{x,z,u\}$ is not a convex geometry, a contradiction.
\end{proof}

\subsection{Proof of the main theorem}

Each trace-minimal non-edge-shelling convex geometry does not satisfy Rule Y. 
The necessity of Theorem \ref{thm:main} follows from the fact that any trace-minimal non-edge-shelling convex geometry is not an edge-shelling convex geometry.

We prove the sufficiency of Theorem \ref{thm:main}. So assume that a convex geometry of stem size 2 satisfying Rule Y is given. Fix an element $x\in \ex(E)$.

\begin{lem}\label{lem:tree}
For a convex geometry of stem size 2 satisfying Rule O and $x\in \ex(E)$, the Hasse diagram of the ordering $>_x$ appeared in Lemma \ref{lem:strict} is a tree except isolated vertices.
\end{lem}

\begin{proof}
Suppose that the Hasse diagram has a cycle. Take a cycle of minimum length. Assign the direction to each edge on the cycle between covering pair of vertices from the smaller vertex to the larger vertex with respect to $>$. Since $>$ is a strict partial order, there exists a source. Note that $u$
is not covered by $x$. Let $y$ and $z$ be vertices covering $u$. By the assumption of the minimality of the cycle, neither $y>_x z$ nor $z>_x y$ holds. Hence we have $(\{x,z\},y)\notin {\cal C}$ and $(\{x,y\},z)\notin {\cal C}$. By applying Rule O to $x,y,z,u$, one of $(\{x,z\},y)\in {\cal C}$ and $(\{x,y\},z)\in {\cal C}$ holds, a contradiction.
\end{proof}

Some element in $E$ may be contained in no rooted circuits containing $x$.

\begin{lem}\label{lem:parallel}
Consider a convex geometry satisfying Rule Y.
Let $y\in E-\{x\}$ be an element such that no rooted circuit contains both $x$ and $y$. Then $(\{x,z\},u)\in {\cal C}$ if and only if $(\{y,z\},u)\in {\cal C}$.
\end{lem}

\begin{proof}
By applying Rule Y to $y$ and $(\{x,z\},u)\in {\cal C}$, either $(\{y,z\},u)\in {\cal C}$ or $(\{y,x\},u)\in {\cal C}$ holds. Since $y$ is not contained in any rooted circuit containing $x$, we have $(\{y,z\},u)\in {\cal C}$.

The converse direction follows similarly.
\end{proof}

We call such an element $y$ an {\it element parallel to} $x$.
Denote the set of elements in $E$ except $x$ and the elements parallel to $x$ by $E_x$ .

Although the vertices in the Hasse diagram correspond to $E_x \cup \{x\}$, we want to make a new tree whose edges correspond to $E$ by using the Hasse diagram.

The Hasse diagram has the maximum element $x$. Every element $y \in E_x$ is covered by a unique element $y'$ by Lemma \ref{lem:tree}. So we can associate every element $y \in E_x$ to the edge $\{y,y'\}$. Add an edge corresponding to $x$ at the maximum element of the graph. Moreover add edges corresponding to the parallel elements to $x$ at vertex $x$.
Then any element in $E$ corresponds to some edge in the graph.

\begin{figure}[ht]
\begin{center}
\includegraphics[bb=0 0 60mm 40mm]{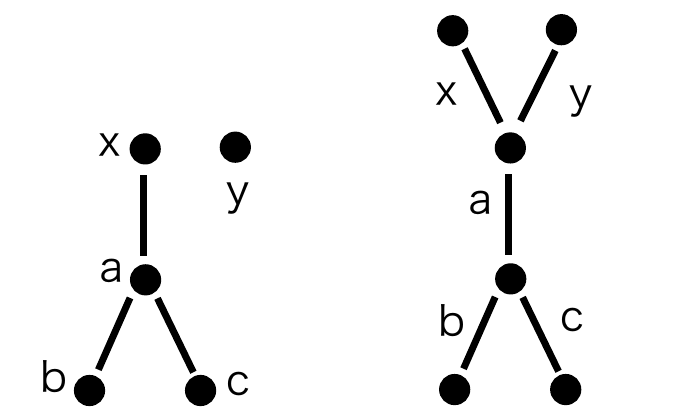}
\caption{Hasse diagram and the corresponding tree}\label{fig:edgeshelling}
\end{center}
\end{figure}

For example, we consider a convex geometry satisfying Rule Y with $(\{x,b\},a), (\{x,c\},a)\in {\cal C}$ as shown in Figure \ref{fig:edgeshelling}(left). The ordering $>_x$, appeared in Lemma \ref{lem:strict}, satisfies $x>_x a>_x b$ and $a>_x c$. Moreover assume that it has $y \in E$ parallel to $x$.  By Lemma \ref{lem:parallel}, it must have $(\{y,b\},a), (\{y,c\},a)\in {\cal C}$. Then we have the corresponding tree as shown in Figure \ref{fig:edgeshelling}(right).

We show that this tree gives rise to the given convex geometry.
For that purpose, consider whether $(\{y,z\},u)$ is a rooted circuit or not for three edges $y,z,u \in E_x$.

Let $u\in E_x$ satisfy that the deletion of $u$ from the tree has two connected components. Then $y$ does not belong to the component containing $x$ if and only if $u>y$. So by definition, it is equivalent to $(\{x,y\},u)\in {\cal C}$. So $y\in E-\{x,u\}$ belongs to the component containing $x$ if and only if $(\{x,y\},u)\notin {\cal C}$.

We say that $y$ belongs to the opposite component to $x$ in this case. 

\begin{lem}\label{lem:oo}
When both $y$ and $z$ belong to the opposite component to $x$, $(\{y,z\},u)$ is not a rooted circuit.
\end{lem}

\begin{proof}
In this case, $(\{x,z\},u)\in {\cal C}$, $(\{x,y\},u)\in {\cal C}$ hold. 
If $(\{y,z\},u)\in {\cal C}$, it contradicts Rule Y to $x$ and $(\{y,z\},u)\in {\cal C}$.
\end{proof}

\begin{lem}\label{lem:ss}
When both $y$ and $z$ belong to the same component to $x$, $(\{y,z\},u)$ is not a rooted circuit.
\end{lem}

\begin{proof}
In this case, $(\{x,y\},u)\notin {\cal C}$ and $(\{x,z\},u)\notin {\cal C}$ hold.
Suppose that $(\{y,z\},u)\in {\cal C}$.
By applying Rule Y to $x$ and $(\{y,z\},u)\in {\cal C}$, either $(\{z,x\},u)\in {\cal C}$ or $(\{y,x\},u)\in {\cal C}$ holds, a contradiction.
\end{proof}

\begin{lem}\label{lem:os}
When $y$ belongs the same component to $x$, and $z$ belongs to the opposite component, $(\{y,z\},u)\in {\cal C}$ holds.
\end{lem}

\begin{proof}
Assume that $(\{x,z\},u)\in {\cal C}$ and $(\{x,y\},u)\notin {\cal C}$.
By applying Rule Y to $y$ and $(\{x,z\},u)\in {\cal C}$, we have $(\{y,z\},u)\in {\cal C}$.
\end{proof}

By Lemmas \ref{lem:oo}, \ref{lem:ss}, \ref{lem:os}, and \ref{lem:parallel}, the tree gives rise to the given convex geometry.
We have completed the proof of the sufficiency.


\section{Other classes of convex geometries of stem size 2}\label{sec:simplicial}

There are many classes of the convex geometries of stem size 2 besides the edge-shelling convex geometry of a tree,
 for example, a double-shelling convex geometry, the vertex-shelling convex geometry of a tree, the simplicial-shelling convex geometry of a chordal graph and so on. 

A simplicial-shelling convex geometry is the convex geometry arising from a chordal graph.
A graph is {\it chordal} if each of its cycles of four or more nodes has a chord, which is an edge joining two nodes that are not adjacent in the cycle. A vertex is called {\it simplicial} if its neighbors are pairwise 
adjacent. The extreme set $\ex(E)$ of the simplicial-shelling convex geometry of a chordal graph consists of all the simplicial vertices of the graph. It is known that the simplicial-shelling convex geometry of a chordal graph is a convex geometry.

The class of simplicial-shelling convex geometries is not closed under trace operation. 

\begin{exa}
We show an example of a simplicial-shelling convex geometry whose trace is not a simplicial-shelling convex geometry. 
Consider the simplicial-shelling convex geometry induced from the chordal graph shown in Figure \ref{fig:chordal}(left).

\begin{figure}[ht]
\begin{center}
\includegraphics[bb=0 0 80mm 45mm]{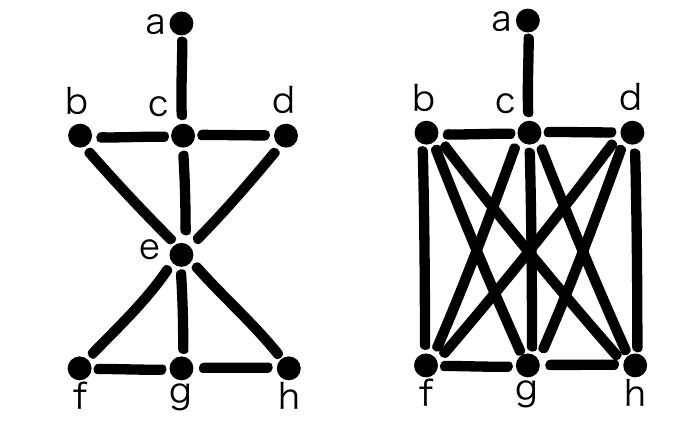}
\caption{A simplicial shelling whose trace is not a simplicial shelling}\label{fig:chordal}
\end{center}
\end{figure}

Assume that the chordal graph gives rise to a simplicial-shelling convex geometry with a family $\cal C$ of rooted circuits.
Consider the trace ${\cal C}:T$ of $\cal C$ on $T=E-\{e\}$. 
Suppose that ${\cal C}:T$ is the rooted circuits of the simplicial-shelling convex geometry of some graph $G'$. 
Note that ${\cal C}:T$ is $\{(X,e)\in {\cal C}:X \subseteq T, e\in T\}$.
Therefore, we can obtain the edges of $G'$ from ${\cal C}:T$.
For example, since $(\{a,b\},c) \in {\cal C}$, $(\{a,b\},c) \in {\cal C}:T$.
This means that $G'$ has edges $\{a,c\}$ and $\{b,c\}$.
By easy calculation, the graph $G'$ must be the graph shown in Figure \ref{fig:chordal}(right). 
Because $G'$ has a chordless cycle $b \to f \to d \to h \to b$, $G'$ is not a chordal graph.
So ${\cal C}:T$ is not a simplicial-shelling convex geometry.
\end{exa}

So a simplicial-shelling convex geometry cannot be characterized by trace-minimal forbidden minors.

The class of convex geometries arising from the vertex shelling of a tree is not closed under trace operation. This class is a subclass of simplicial-shelling convex geometries.

The class of double-shelling convex geometries is another important class of convex geometries of stem size 2. A double-shelling convex geometry is induced from a poset. This class is closed under trace operation.
Kashiwabara and Nakamura \cite{KN} succeeded in characterizing double-shelling convex geometries in terms of trace-minimal forbidden minors.

\section{Characterization of the convex geometries of stem size 2}\label{sec:eminor}

In this section, we characterize the class of convex geometries of stem size 2 in terms of the trace-minimal forbidden minors.

\begin{df}
A family of rooted circuits of stem size 2 is called a {\it trace-minimal non-convex-geometry} of stem size 2 if it satisfies the following two conditions:
(1) It is the rooted circuits of no convex geometry,
(2) Its trace on any proper subset is the rooted circuits of a convex geometry.
\end{df}

\begin{lem}\label{lem:stemclose}
The class of convex geometries of stem size 2 is closed under trace operation.
As for the class of rooted circuits of stem size 2, any trace-minimal non-convex-geometry has a ground set of at most size 5.
\end{lem}


\begin{proof}
It follows from Lemma \ref{lem:tracerooted} that the class is closed under trace operation. 

Let $\cal C$ be a set of rooted circuits of stem size 2.
Assume that  $\cal C$ is the rooted circuits of no convex geometry.
Moreover, suppose that there exists a trace-minimal non-convex-geometry of size more than 5.
Then, by the axiom of convex geometries in terms of rooted circuits, there exist 5 points $\{x,y,z,u,v\}$ in the trace-minimal non-convex-geometry such that $(\{x,y\},z)$ and $(\{z,v\},u)$ belong to $\cal C$ by Lemma \ref{lem:rootedaxiom}, and 
 neither $(\{x,v\},u)$ nor $(\{y,v\},u)$ belongs to $\cal C$.
The trace of $\cal C$ on $\{x,y,z,u,v\}$ is not a convex geometry, a contradiction.
\end{proof}

By specifying all trace-minimal non-convex geometries of stem size 2, we characterize the class of convex geometries of stem size 2. After that, we may consider the non-convex-geometries of stem size 2 as the trace-minimal forbidden minors.

\begin{figure}[ht]
\begin{center}
\includegraphics[bb=0 0 100mm 30mm]{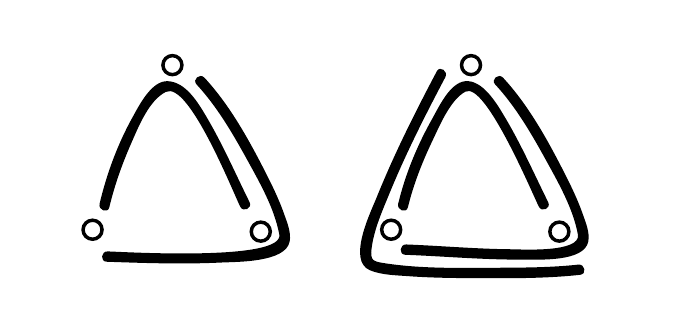}
\caption{Trace-minimal non-convex-geometries of size 3}\label{fig:e3}
\end{center}
\end{figure}

\begin{figure}[ht]
\begin{center}
\includegraphics[bb=0 0 100mm 30mm]{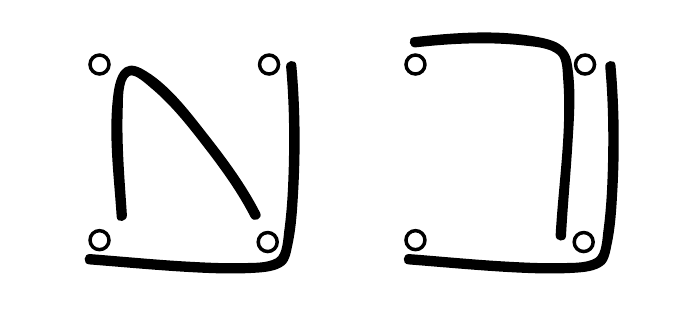}
\caption{Trace-minimal non-convex-geometries of size 4 with two rooted circuits}\label{fig:e4-2}
\end{center}
\end{figure}

\begin{figure}[ht]
\begin{center}
\includegraphics[bb=0 0 150mm 30mm]{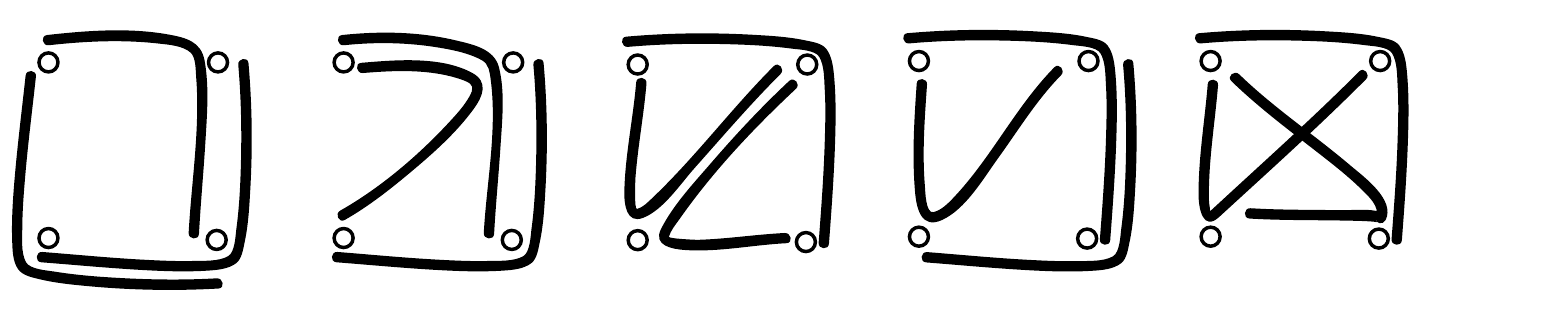}
\caption{Trace-minimal non-convex-geometries of size 4 with three rooted circuits}\label{fig:e4-3}
\end{center}
\end{figure}

\begin{figure}[ht]
\begin{center}
\includegraphics[bb=0 0 150mm 30mm]{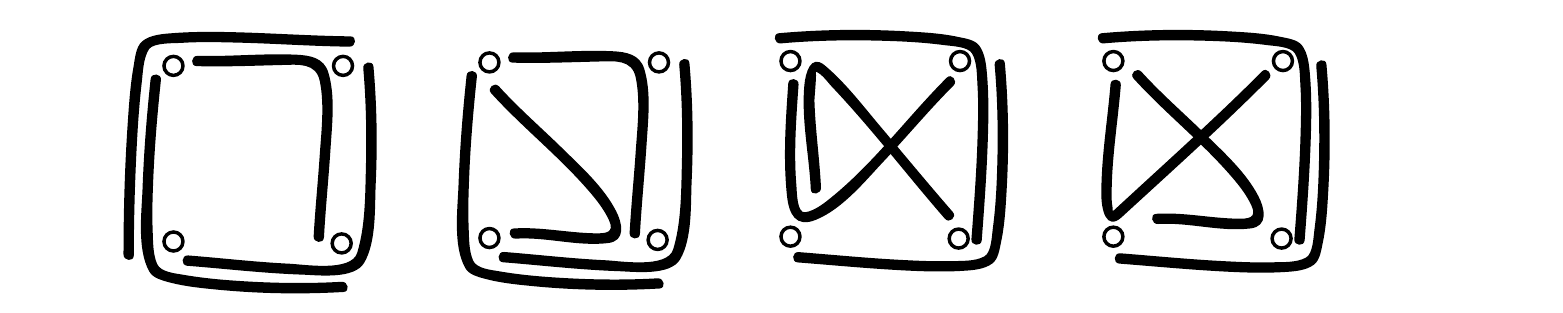}
\caption{Trace-minimal non-convex-geometries of size 4 with four rooted circuits}\label{fig:e4-4}
\end{center}
\end{figure}

\begin{figure}[ht]
\begin{center}
\includegraphics[bb=0 0 65mm 35mm]{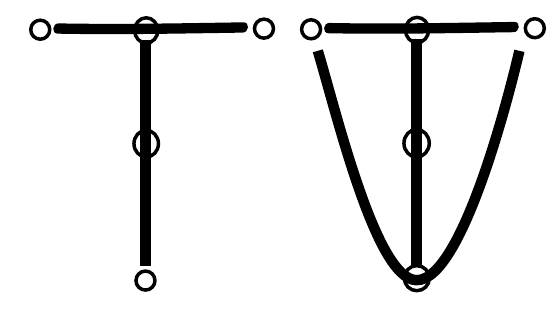}
\caption{Trace-minimal non-convex-geometries of size 5}\label{fig:e5}
\end{center}
\end{figure}

\begin{thm}
The trace-minimal non-convex-geometries for the class of convex geometries of stem size 2 are all identified as listed in Figures \ref{fig:e3}, \ref{fig:e4-2}, \ref{fig:e4-3}, \ref{fig:e4-4}, and  \ref{fig:e5}.
\end{thm}

Any trace-minimal non-convex-geometry of stem size 2 has size at most 5 by Lemma \ref{lem:stemclose}. 

We make the list up to permutation of the ground set. 

First, we make families of rooted circuits of small size. 
We obtain a new family by adding a rooted circuit to them iteratively, identify isomorphic families, and pick up trace-minimal non-convex-geometries.

There exist two of size 3 shown in Figure \ref{fig:e3}.
Note that when there exist more than one rooted circuit in three points, the family contains a trace-minimal non-convex-geometry of size 3.
So since any trace-minimal non-convex-geometry of size at least 4 does not contain these minors, the ground set $\{x,y,z\}$ of  each rooted circuit $(\{x,y\},z)$  is different from one another.
This fact leads us to  identify trace-minimal non-convex-geometries easier.

There are many trace-minimal non-convex-geometries of size 4. So we classify them according to the number of their rooted circuits.

There are two families which have two rooted circuits, shown in Figure \ref{fig:e4-2}.

There are five families which have three rooted circuits, shown in Figure \ref{fig:e4-3}.

There are four families which have four rooted circuits, shown in Figure \ref{fig:e4-4}.

There are two families of size 5, shown in Figure \ref{fig:e5}.

\end{document}